\documentclass[12pt]{article}

\usepackage{amsmath}
\usepackage{amsfonts}
\usepackage[T1]{fontenc} 
\usepackage[dvips,pdftex]{graphicx}
\usepackage{a4wide}

\newtheorem{prop}{Proposition}[section]
\newtheorem{corol}{Corollary}[section]

\newtheorem{lemme}{Lemma}[section]

\newcommand{\CQFD}{\hfill $\square$}

\newcommand{\ind}{\mathbf{1}}

\def\real{\Bbb{R}}

\def\E{\mathop{\hbox{\rm I\kern-0.20em E}}\nolimits}

\def\og{\leavevmode\raise.3ex
     \hbox{$\scriptscriptstyle\langle\!\langle$~}}
\def\fg{\leavevmode\raise.3ex
     \hbox{~$\!\scriptscriptstyle\,\rangle\!\rangle$}~}
     
\begin{document}

\title{Some remarks on Betti numbers of random polygon spaces}
\author{ Clément Dombry\thanks{Laboratoire de Math\'ematiques et Applications, T\'el\'eport 2- BP30179, Boulevard Pierre et Marie Curie, 86962 Futuroscope Chasseneuil Cedex, France. clement.dombry@math.univ-poitiers.fr} and Christian Mazza\thanks{D\'epartement de Math\'ematique, Université de Fribourg, Chemin du Mus\'ee 23, CH-1700 Fribourg, Suisse.Christian.Mazza@unifr.ch} \protect\hspace{1cm}}
\maketitle

\abstract{Polygon spaces like $M_\ell=\{(u_1,\cdots,u_n)\in S^1\times\cdots S^1\ ;\ \sum_{i=1}^n l_iu_i=0\}/SO(2)$
or they three dimensional analogues $N_\ell$ 
play an important rôle in geometry and topology, and are also of  interest in robotics
where the $l_i$ model the lengths of robot arms. When $n$ is large, one can assume that each
$l_i$ is a positive real valued random variable, leading to a random manifold. The complexity
of such manifolds can be approached by computing  Betti numbers, the Euler characteristics, or
the related Poincar\'e polynomial. We study the average values of Betti numbers of dimension $p_n$
when $p_n\to\infty$ as $n\to\infty$. We also focus on the limiting mean Poincar\'e polynomial, in two and
three dimensions. We show that in two dimensions, the mean total Betti number behaves as the total Betti number
associated with the equilateral manifold where $l_i\equiv \bar l$. In three dimensions, these two quantities
are not any more asymptotically equivalent. We also provide asymptotics for the Poincar\'e polynomials}
\medskip


\section{Introduction}
\subsection{Background}
In this note, we consider a question raised by M.Farber in \cite{Fa}.
We study closed planar $n$-gons whose sides have fixed lengths $l_1,\cdots,l_n$ where $l_i>0$ for $1\leq i\leq n$. The set of polygonal linkage in ${\mathbb R}^2$
$$
M_\ell=\{(u_1,\cdots,u_n)\in S^1\times\cdots S^1\ ;\ \sum_{i=1}^n l_iu_i=0\}/SO(2)
$$
parametrizes the variety of all possible shapes of such planar $n$-gons with sides given by $\ell=(l_1,\cdots,l_n)$. The unit vector $u_i\in{\mathbb C}$ indicates the direction of the $i$-th side of the polygon. The condition $\sum l_iu_i=0$ expresses the property of the polygon being closed. The rotation group $SO(2)$ acts on the set of side directions $(u_1,\cdots,u_n)$ diagonnaly.\\

Polygon spaces
play a fundamental rôle in topology and geometry, as illustrated for example 
by Kempe Theorem which states that "{\it Toute courbe alg\'ebrique peut \^etre trac\'ee \`a l'aide d'un
syst\`eme articul\'e}", see e.g. \cite{Lebe}. 
\cite{Fa3} provides other examples of such universality results in 
topology. Polygon spaces generated an active research area in geometry
(see e.g. \cite{FaSc}, \cite{Hausmann}, or \cite{Thurston}), but are also
of strong
interest in applications like robotics where each $l_i$ models the length or a robot arm (see e.g. \cite{Fa}, \cite{FaKa}, \cite{FaSc} and \cite{Fa3}).
We can 
also point out potential applications in polymer science where such polygons model proteins. In systems
composed of a large number $n\!>\!\!>\!1$ of components, the $l_i$ are usually only
partially known, so that we can assume that each $l_i\in\real^+$ is random.
We denote by $\mu_n$ the 
distribution of $\ell$. We will obtain our results under the following assumption:
\begin{center}
({\bf H}) $\mu_n$ is a product measure $\mu_n=\mu^{\otimes n}$ with $\mu$ a diffuse measure on $(0,\infty)$ such that
$$\int e^{\eta  x}\mu(dx)<\infty {\rm\ \ for\  some\ \ } \eta >0.$$
\end{center}
Notice that $M_{t\ell}$ and
$M_\ell$ are equal when $t>0$, so that the measure $\mu_n$ might be seen as a probability
measure on the unit simplex $\triangle^{n-1}$.

To get some idea on the nature of the random manifold $M_\ell$, one can study the stochastic behavior of invariants
like Betti numbers, the Euler characteristics or the total Betti number (see below).
Here, we focus on the  Betti numbers $b_p(M_\ell)$,
for dimensions $p=p_n$ growing with $n$.
 We recall the result of \cite{FaSc} describing Betti numbers of planar polygon spaces as functions of the length vector $\ell$. In what follows,
$[n]$ denotes the set $\{1,\cdots,n\}$.
A subset $J\subset [n]$ is called {\it short} if
$$
\sum_{j\in J} l_j<\sum_{j\notin J} l_j.
$$
It is called {\it median} if 
$
\sum_{j\in J} l_j=\sum_{j\notin J} l_j
$.
Let $1\leq i_0\leq n$ be such that $l_{i_0}$ is maximal among $l_1,\cdots,l_n$. Denote by $a_p(\ell)$ the number of short subsets $J\subset [n]$ of cardinality $|J|=p+1$ and containing $i_0$. Denote by $\tilde a_p(\ell)$ the number of median subsets $J\subset [n]$ of cardinality $|J|=p+1$ and containing $i_0$. Then one has for $p=0,1,\dots, n-3$
\begin{equation}\label{eq:def-Betti}
b_p(M_\ell)=a_p(\ell)+\tilde a_p(\ell)+a_{n-3-p}(\ell),
\end{equation}
so that the Poincar\'e polynomial of the random manifold is given by
\begin{equation}\label{PC2}
p_{M_\ell}(t) = \sum_{p=0}^{n-3}b_p(M_\ell)t^p
 = q(t)+t^{n-3}q(\frac{1}{t})+r(t),
\end{equation}
where
$$q(t)=\sum_{k=0}^{n-3}a_k t^k\hbox{ and }r(t)=\sum_{k=0}^{n-3}\tilde a_k t^k,$$
see \cite{Fa}. The total Betti number $B(M_\ell)$ defined by
$$B(M_\ell)=\sum_{p=0}^{n-3}b_p(M_\ell)=p_{M_\ell}(1),$$
provides ideas on the size or on the complexity of the manifold $M_\ell$. We will study the
asymptotic behavior of $B(M_\ell)$ when $n$ is large and $\ell$ is random. We first give some examples
following \cite{FaSc}. 
\medskip

In the equilateral case where each $l_i$ is equal to some $\bar l>0$, it turns out that one can give
exact formulas for the various Betti numbers, and therefore for $B(M_{\bar\ell})$: assume $n=2r+1$ odd. Then
$b_k(M_{\bar \ell})={n-1\choose k}$ when $k<r-1$, $b_k(M_{\bar \ell})=2{n-1\choose r-1}$ when $k=r-1$ and
$b_k(M_{\bar \ell})={n-1\choose k+2}$ when $k>r-1$. The related total Betti number is then given by
$p_{M_{\bar \ell}}(1)=2^{n-1}-{n-1\choose r}$. For arbitrary large $n$, one has (see \cite{FaSc})
\begin{equation}\label{AsymptoticEquilateral}
p_{M_{\bar \ell}}(1) = B_n = 2^{n-1}(1-\sqrt{\frac{2}{\pi n }}+o(n^{-1/2}))\ \ ,\ \ n\to\infty.
\end{equation}
For pentagons, that is when $n=5$, the moduli space $M_\ell$ is a compact orientable
surface of genus not exceeding 4. 
\medskip

The vector length $\ell$ is said to be {\it generic} when $\sum_{i=1}^n l_i\varepsilon_i\ne 0$, for any
$\varepsilon =(\varepsilon_i)_{1\le i\le n}$, where , $\forall i\ \varepsilon_i\in \{-1,+1\}$.
When $n$ is even, equilateral weights with $l_i\equiv \bar l$ are not generic. \cite{FaSc}
proved that for generic $\ell$, the total Betti number $B(M_\ell)$ is bounded by $2B_{n-1}$, so that
the explicit formulas obtained for
equilateral $n$-gons provide bounds for the maximum over $\ell$ of $B(M_\ell)$.
\bigskip

\subsection{Results\label{s.results}}
\bigskip

\cite{FaKa} considered the special case where $\mu$ is the uniform probability measure on the unit interval
with $\mu_n=\mu^{\otimes n}$, and the case where $\mu_n$ is the uniform measure on the simplex $\Delta^{n-1}$.
It was proven that for fixed $p\geq 0$,
the average $p$-dimensional betti number 
$$
\mu_n[b_p(M_\ell)]=\int b_p(M_\ell)\mu_n(d\ell)
$$ 
is asymptotically equivalent to ${n-1}\choose{p}$, the difference going to zero at an exponential speed. The techniques use exact formulas
for the volume of the intersection of a half space with a simplex.
We will avoid such formulas to treat general diffuse probability measures using
probabilistic techniques, since in fact such volume formulas do not exist for arbitrary measures. 
Next, \cite{Fa} consider both planar and spatial polygon
spaces, and proved, under an admissibility condition on $\mu_n$ 
similar results for mean Betti numbers and also for higher moments, again for fixed dimensions $p$. 
As an open question, the author raises the issue of computing the average total Betti number
$$
\mu_n[B(M_\ell)]=\int B(M_\ell)\mu_n(d\ell).
$$
We will consider more generally the mean Poincar\'e polynomial
$$\bar p_{M_\ell}(t)=\mu_n[p_{M_\ell}(t)] = \int p_{M_\ell}(t)\mu_n(d\ell),$$
with $\bar p_{M_\ell}(1)=\mu_n[B(M_\ell)]$.
As the author notices, the knowledge of the individual average Betti numbers $\mu_n[b_p(M_\ell)]$ for large $n$ and fixed $p$ can't help since the terms cannot simply be added up. We will therefore  consider the asymptotic behavior
of high dimensional Betti numbers $\mu_n[b_{p_n}(M_\ell)]$, where $p_n$ goes to infinity when $n\to\infty$ (see Proposition \ref{prop:pr1}).

We will obtain our results for product measure $\mu_n$ satisfying assumption $({\bf H)}$ and assume throughout the paper that this hypothesis is satisfied.
We  prove in Proposition \ref{prop:pr2} that  the mean total Betti number is such that
\begin{equation}\label{FirstAsymp}
\bar p_{M_\ell}(1) = \mu_n[B(M_\ell)]\sim 2^{n-1},
\end{equation}
This shows that equilateral polygons (see (\ref{AsymptoticEquilateral})) are representative
of the emerging average manifold as $n>>1$, as suggested in \cite{Fa}.
We will also consider the mean Poincar\'e polynomial
as $n$ is large, and show that 
$$\bar p_{M_\ell}(t)\sim (1+t)^{n-1}\hbox{ when } 0<t<1,$$
and that
$$\bar p_{M_\ell}(t)\sim (1+t)^{n-1}t^{-2}\hbox{ when } t>1.$$
Further moments are also considered and their asymptotic is given in Proposition \ref{prop:pr4}.

We next consider spatial polygon spaces
$$
N_\ell=\{(u_1,\cdots,u_n)\in S^2\times\cdots S^2\ ;\ \sum_{i=1}^n l_iu_i=0\}/SO(3).
$$
In this case, for generic length vector $\ell$, \cite{HausmannKnu} proved that the even Betti numbers are given by
\begin{equation}\label{BettiSpatial}
b_{2p}(N_\ell)=\sum_{j=0}^p (\hat a_j(\ell)-\hat a_{n-j-2}(\ell)),
\end{equation}
where $\hat a_j(\ell)$ denotes the number of short subsets $J\subset [n]$ of cardinality $|J|=j+1$ containing $n$.
The Betti number of odd dimensions vanish. 
Furthermore, \cite{HausmannKnu} proved that the related Poincar\'e polynomial is given by
\begin{equation}\nonumber
 p_{N_\ell}(t)=\frac{1}{1-t^2}\Big(\sum_{J\in {\cal S}_n}t^{2(\vert J\vert-1)}-t^{2(n-\vert J\vert -1)}\Big),
\end{equation}
where $J\in{\cal S}_n$ if and only if $\{n\}\subset J\subset [n]$  and $J$ is median or short. If $\ell$ is generic, there is no median set and this is equivalent to 
\begin{equation}\label{PC3}
 p_{N_\ell}(t)=\frac{1}{1-t^2}\sum_{j=0}^{n-1}\hat a_j(t^{2j}-t^{2(n-j-2)})=\frac{1}{1-t^2}\left[\hat q(t^2)-t^{2(n-2)}\hat  q(t^{-2}) \right],
\end{equation}
where
$$\hat q(t)=\sum_{j=0}^{n-1} \hat a_j t^j.$$
\medskip

In the equilateral case where $l_i\equiv \bar l$, \cite{Kl} proved that 
the $2p$-dimensional Betti number $b_{2p}(N_{\bar\ell})$ is given by
$$b_{2p}(N_{\bar\ell})=\sum_{i=0}^p {n-1\choose i},$$
when $n=2k+1\ge 3$, so that
the Euler characteristics or total Betti number is
explicitely given as
$$p_{N_{\bar\ell}}(1)=\sum_{i=0 }^{k-1}{2k \choose i}(k-i),$$
with
$$p_{N_{\bar\ell}}(1)\sim \sqrt{\frac{n}{2\pi}}2^{n-2}.$$

We will study the asymptotic behavior of the mean Poincar\'e polynomial 
$$\bar p_{N_\ell}(t)=\mu_n[p_{N_\ell}(t)] = \int p_{N_\ell}(t)\mu_n(d\ell),$$
in the large $n$
limit by providing large deviations estimates. We will see that
$$\bar p_{N_\ell}(1) = \mu_n[B(N_\ell)] \sim n 2^{n-2} >\!\!> p_{N_{\bar\ell}}(1).$$
Furthermore, we will see
in Proposition \ref{prop:pr3} that
the mean Poincar\'e polynomial exhibits asymptotically a singular behavior in the neighborhood of $t=1$, that is
$$\bar p_{N_\ell}(t)\sim \frac{(1+t^2)^{n-1}}{(1-t^2)}\hbox{ when } 0<t<1,$$
and
$$\bar p_{N_\ell}(t)\sim \frac{(1+t^2)^{n-1}}{(t^2-1)t^2}\hbox{ when } t>1.$$
This shows that equilateral configuration spaces are not representative of the random manifold
in dimension 3 when $n$ is large.

\section{Preliminaries}
We introduce here the main technical tool used in our analysis of the Betti numbers of random polygon spaces: a probabilistic interpretation of formulas (\ref{eq:def-Betti}) and (\ref{BettiSpatial}) in terms of random permutations and  stopping times.
We first introduce some notations.\\

\noindent
For any length vector $\ell\in (0,\infty)^n$, we define
$\tilde{\ell}$ obtained from $\ell$ by the following permutation
 of the coordinates : let $i_0$ be the minimal index such that $l_{i_0}$ is maximal among the $l_i, 1\leq i\leq n$, and define $\tilde \ell=(\tilde l_1,\cdots,\tilde l_1)$ by $\tilde l_n=l_{i_0}$, $\tilde l_{i_0}=l_{n}$, and $\tilde l_{i}=l_i$ if $i\notin\{i_0,n\}$.\\
We denote by $\sigma$ a random permutation of $\Sigma_{n-1}$ with uniform distribution ${\cal U}_{\Sigma_{n-1}}$.
The stopping time $\tau_\sigma(\ell)$ is defined by
$$
\tau_\sigma(\ell)=\min\left\{0\leq t\leq n-1\ ;\  \sum_{i=1}^t l_{\sigma(i)}+ l_n-\sum_{i=t+1}^{n-1} l_{\sigma(i)}\geq 0\right\}.
$$
We use also the notation $\tau(\ell)=\tau_{Id}(\ell)$ and $\tilde\tau(\ell)=\tau_{Id}(\tilde\ell)$. Please note that these stopping times are well-defined and that $\tau\leq n-1$ and $\tilde\tau\leq n-2$. \\
We denote by $k$ a random variable with binomial distribution ${\cal B}_{n-1,q}$ with parameters $n-1$ and $q\in[0,1]$.\\

First consider the planar case.
\begin{lemme}\label{lem:representation1}
The number $a_p(\ell)$ of short sets is given by
$$ 
a_p(\ell)={{n-1}\choose p}{\cal U}_{\Sigma_{n-1}}[\tau_\sigma(\tilde\ell)>p]
$$
The number of median sets $\tilde a_p(\ell)$ vanishes $\mu_n$-almost surely.\\ 
Hence, the planar Betti numbers are given $\mu_n$-almost surely by
$$
b_p(M_\ell)={\cal U}_{\Sigma_{n-1}}\left[{{n-1}\choose p}\ind_{\{\tau_\sigma(\tilde\ell)>p\}}+{{n-1}\choose {p+2}}\ind_{\{\tau_\sigma(\tilde\ell)>n-p-3\}}\right]
$$ 
and have expected value
$$
\mu_n[b_p(M_\ell)]=\mu_n\left[{{n-1}\choose p}\ind_{\{\tilde\tau(\ell)>p\}}+{{n-1}\choose {p+2}}\ind_{\{\tilde\tau(\ell)>n-p-3\}}\right]
$$
\end{lemme}

In the spatial case, the following representation holds.
\begin{lemme}\label{lem:representation2}
The coefficients $\hat a_p$ are given by
$$ 
\hat a_p(\ell)={{n-1}\choose p}{\cal U}_{\Sigma_{n-1}}[\tau_\sigma(\ell)>p]
$$
Hence, the even spatial Betti numbers are given $\mu_n$-almost surely by
$$
b_{2p}(N_\ell)=2^{n-1}({\cal U}_{\Sigma_{n-1}}\otimes B_{n-1,1/2})\left[\ind_{\{\tau_\sigma(\ell)>k;0\leq k\leq p \}}- \ind_{\{\tau_\sigma(\ell)>k;n-p-2\leq k\leq n-2 \}}\right]
$$ 
and have expected value
$$
\mu_n[b_{2p}(N_\ell)]=2^{n-1}(\mu_n \otimes B_{n-1,1/2})\left[\ind_{\{\tau(\ell)>k;0\leq k\leq p \}}- \ind_{\{\tau(\ell)>k;n-p-2\leq k\leq n-2 \}}\right]
$$
\end{lemme}

\noindent
{\bf Proof of Lemmas \ref{lem:representation1} and \ref{lem:representation2}}\\ 
We consider the planar case and prove the first lemma. The second lemma corresponding to the spatial case is proved with a very similar analysis.\\
\noindent From the definition of the coefficient $a_p(\ell)$,  we have
$$a_p(\ell)=\sum_{J\subset [n-1]\ ;\ |J|=p} \ind_{A_J}(\tilde \ell)$$ 
where $A_J=\{ \ell \ ;\ \sum_{j\in J} l_j+ l_n-\sum_{j\notin J} l_j<0\}$. From the definition of $\tau_\sigma$, it is easily seen that  $\ell\in A_{\{\sigma(1),\cdots,\sigma(p)\}}$ if and only if $\tau_\sigma(\ell)>p$. Furthermore, for each subset $J\subset [n-1]$ such that $|J|=p$, there are $p!(n-1-p)!$ permutations $\sigma\in\Sigma_{n-1}$ such that $J=\{\sigma(1),\cdots,\sigma(p)\}$. As a consequence, the coefficient $a_p(\ell)$ rewrites
\begin{eqnarray*}
a_p(\ell)&=&\frac{1}{p!(n-1-p)!}\sum_{\sigma\in\Sigma_{n-1}} \ind_{\{ \tau_\sigma(\tilde\ell)> p\}}\\
&=& {{n-1}\choose p}{\cal U}_{\Sigma_{n-1}}[\tau_\sigma(\tilde\ell)>p].
\end{eqnarray*}
From the definition of the coefficient $\tilde a_p(\ell)$,  we have
$$\tilde a_p(\ell)=\sum_{J\subset [n-1]\ ;\ |J|=p} \ind_{B_J}(\tilde \ell)$$ 
where $B_J=\{ \ell \ ;\ \sum_{j\in J} l_j+ l_n-\sum_{j\notin J} l_j=0\}$. But it is easily seen that since $\mu$ is diffuse, the sum $\sum_{j\in J} l_j+ l_n-\sum_{j\notin J} l_j$ is also diffuse
and $\mu_n(B_J)=0$ for any $J\subset [n-1]$. Hence $\tilde a_p(\ell)$ is almost surely equal to zero.\\
 
\noindent The formula for the Betti number $b_p(\ell)$ is then a reformulation of equation (\ref{eq:def-Betti}). Thanks to the invariance of $\mu_n$ under the action of the permutation group, the distribution of $\tau_\sigma(\tilde \ell)$ under $\mu_n$ does not depend on $\sigma\in\Sigma_{n-1}$ and hence is equal to the distribution of $\tilde\tau$. The result for the expected value $\mu_n[b_p(M_\ell)]$ follows.
\CQFD \\
\ \\

As will be clear in the sequel, the asymptotic behavior of the Betti numbers is strongly linked with the asymptotic behavior of the random variables $\tau(\ell)$ and $\tilde\tau(\ell)$. This is the point of the following lemma.
\begin{lemme}\label{lem:LLN-TCL}
The following weak convergences hold under $\mu_n$ as $n\to\infty$:
\begin{enumerate}
\item weak law of large numbers:
$$ n^{-1}\tau \Rightarrow \delta_{1/2}, $$
\item central limit theorem:
$$n^{-1/2}\left(\tau-\frac{n}{2}\right) \Rightarrow {\mathcal N}(0,\sigma_{\tau}^2),$$
where $\sigma_{\tau}=\frac{\sigma}{2m}$, $m=\E(l)$ and $\sigma^2= {\rm Var}(l)$.
\item large deviations: for any $\varepsilon>0$,
$$\limsup n^{-1}\log \mu_n(\vert n^{-1}\tau-1/2\vert \geq \varepsilon) <0 $$
\end{enumerate}
The same results also hold for $\tilde\tau$ instead of $\tau$ with the same variance $\sigma_{\tilde\tau}=\sigma_\tau$.
\end{lemme}
The proof is postponed to the appendix.

\section{High dimensional Betti numbers\label{HB}}
\subsection{Planar polygons}
The following proposition gives the asymptotic of  average high dimensional Betti numbers.

\begin{prop}\label{prop:pr1}
 Let $(p_n)_{n\geq 1}$ be a sequence of integers.
\begin{enumerate}
\item If $\limsup n^{-1}p_n<1/2$, then $\mu_n[b_{p_n}(M_\ell)]\sim {{n-1}\choose{p_n}}$ as $n\to\infty$.
\item If $\liminf n^{-1}p_n>1/2$, then $\mu_n[b_{p_n}(M_\ell)]\sim {{n-1}\choose{p_n+2}}$ as $n\to\infty$.
\item If $\lim n^{-1/2}(p_n-n/2)=\alpha$, then $\mu_n[b_{p_n}(M_\ell)]\sim \sqrt{\frac{2}{\pi n}}e^{-2\alpha^2} 2^{n-1}$ as $n\to\infty$.
\end{enumerate}
\end{prop}

Applying Proposition \ref{prop:pr1} with a specific choice of the sequence $p_n$, we deduce the following corollary.
The asymptotic of the binomial coefficient is obtained with Stirling's formula.
\begin{corol}
Let $p\in (0,1)$ and $p_n=[np]$. Then, 
$$\lim_{n\to\infty} n^{-1}\log \mu_n[b_{p_n}(M_\ell)] = -p\log p - (1-p)\log(1-p)$$
\end{corol}

\noindent {\bf Proof of Proposition \ref{prop:pr1}:}\\
From Lemma \ref{lem:representation1}, the average Bett numbers is given by
$$\mu_n[b_{p_n}(M_\ell)]={{n-1}\choose{p_n}} \mu_n(\tilde\tau>p_n)+{{n-1}\choose{p_n+2}} \mu_n(\tilde\tau>n-3-p_n).$$
When $\limsup n^{-1}p_n<1/2$, the weak law of large numbers provided in Lemma \ref{lem:LLN-TCL} implies that
$$\mu_n(\tilde\tau>p_n)\rightarrow 1 \ \ {\rm and}\ \ \mu_n(\tilde\tau>n-3-p_n)\rightarrow 0$$
as $n\to\infty$, and from  large deviations estimates, the convergence speed to zero is exponential. The first point in Proposition \ref{prop:pr1} follows since
\begin{eqnarray*}
& &b_{p_n}(n,\mu_n)\\
&=&{{n-1}\choose{p_n}}\left( \mu_n(\tilde\tau>p_n)+\frac{(n-p_n-1)(n-p_n-2)}{(p_n+1)(p_n+2)} \mu_n(\tilde\tau>n-3-p_n)\right)\\
&\sim &{{n-1}\choose{p_n}}.
\end{eqnarray*}

\noindent Similarly, when $\liminf n^{-1}p_n>1/2$, 
$$\mu_n(\tilde\tau>p_n)\rightarrow 0 \ \ {\rm and}\ \ \mu_n(\tilde\tau>n-3-p_n)\rightarrow 1$$
as $n\to\infty$ where the convergence to zero is exponentially fast. The second point in Proposition \ref{prop:pr1} follows.\\

\noindent
Finally, in the case $\lim n^{1/2}(p_n-n/2)=\alpha$, the central limit Theorem from Lemma \ref{lem:LLN-TCL} yields that as $n\to\infty$
$$\mu_n(\tilde\tau>p_n)\rightarrow 1-F_{{\mathcal N}}(\alpha/\sigma_{\tilde\tau}) \ \ {\rm and}\ \  \mu_n(\tilde\tau>n-3-p_n)\rightarrow1-F_{{\mathcal N}}(-\alpha/\sigma_{\tilde\tau}),$$
where $F_{{\mathcal N}}$ is the repartition function of the standard normal distribution.
Furthermore, from the local limit theorem for the binomial distribution, 
$${{n-1}\choose{p_n}}\sim  {{n-1}\choose{n-3-p_n}} \sim \sqrt{\frac{2}{\pi n}}e^{-2\alpha^2} 2^{n-1}, $$
as $n\to\infty$. 
These estimates yield the last point in Proposition \ref{prop:pr1} since 
$$1-F_{{\mathcal N}}(\alpha/\sigma_{\tilde\tau})+1-F_{{\mathcal N}}(-\alpha/\sigma_{\tilde\tau})=1.$$
\CQFD

\subsection{Spatial polygons\label{Admissible}}
We perform a similar study in the spatial case. The asymptotic behavior of average high dimensional Betti numbers is given by the following Proposition.

\begin{prop}\label{prop:pr1-spatial}
Let $(p_n)_{n\geq 1}$ be a sequence of integers. 
\begin{enumerate}
\item If $\limsup n^{-1}p_n<1/2$, then $\mu_n\left[b_{2p_n}(N_\ell)\right]\sim \sum_{k=0}^{p_n}{n-1\choose{k}}$ as $n\to\infty$.
\item If $\liminf n^{-1}p_n>1/2$, then $\mu_n\left[b_{2p_n}(N_\ell)\right]\sim \sum_{k=0}^{n-p_n-3}{n-1\choose{k}}$ as $n\to\infty$.
\item If $\lim n^{-1/2}(p_n-n/2)=\alpha$, then $\mu_n\left[b_{2p_n}(N_\ell)\right]\sim C(\alpha) 2^{n-1}$ as $n\to\infty$, with
$$C(\alpha)=\int_{2|\alpha|}^\infty \frac{e^{-\frac{u^2}{2}}}{\sqrt{2\pi}}\mathbb P(\vert Z\vert <\frac{um}{\sigma}){\rm d}u, $$
where $m=\mu(l)$, $\sigma^2 = {\rm Var}(l)$, and $Z$ is standard normal.
\end{enumerate}
\end{prop}

\noindent {\bf Proof of Proposition \ref{prop:pr1-spatial}:}\\
Recall from Lemma \ref{lem:representation2} that the expected Betti number $\mu_n[b_{2p}(N_\ell)]$ is given by 
$$ \mu_n[b_{2p}(N_\ell)]= 2^{n-1} (\mu_n\otimes{\cal B}_{n-1,1/2})\left[{\bf 1}_{\{\tau >k ; 0\leq k\leq p_n \}}- {\bf 1}_{\{\tau >n-1-k ; 1\leq k\leq p_n+1 \}}\right]$$
(we use here the fact that $k$ and $n-1-k$ have the same distribution under ${\cal B}_{n-1,1/2}$).\\
Consider first the case $\limsup n^{-1}p_n<1/2$ and write 
$$\mu_n(\tau>p_n){\cal B}_{n-1,1/2}(0\leq k\leq p_n)\leq (\mu_n\otimes{\cal B}_{n-1,1/2})\left[{\bf 1}_{\{\tau >k ; 0\leq k\leq p_n \}}\right]\leq {\cal B}_{n-1,1/2}(0\leq k\leq p_n) .$$
Using the weak law of large numbers $n^{-1}\tau\to 1/2$ under $\mu_n$ and the asymptotic for $p_n$, wee see that $\mu_n(\tau>p_n)\to 1$ as $n\to\infty$.  Hence the equivalent 
$$(\mu_n\otimes{\cal B}_{n-1,1/2})\left[{\bf 1}_{\{\tau >k ;0\leq k\leq p_n \}}\right]\sim {\cal B}_{n-1,1/2}(0\leq k\leq p_n).$$
In the same way, 
$$0\leq (\mu_n\otimes{\cal B}_{n-1,1/2})\left[{\bf 1}_{\{\tau >n-1-k ; 1\leq k\leq p_n+1 \}}\right]\leq \mu_n(\tilde\tau>n-1-p_n){\cal B}_{n-1,1/2}(1\leq k\leq p_n+1))$$
and a large deviations argument shows that $\mu_n(\tau>n-p_n)$ converges exponentially fast to zero, so that this last term is of smaller order than ${\cal B}_{n-1,1/2}(0\leq k\leq p_n)$. This proves the first point.\\

\noindent
Consider now the case $\liminf n^{-1}p_n>1/2$. It appears that many terms cancel out and we have for large $n$
\begin{eqnarray*}
\mu_n\left[b_{2p_n}(N_\ell)\right]&=& 2^{n-1} (\mu_n\otimes{\cal B}_{n-1,1/2})\left[{\bf 1}_{\{\tau >k ; 0\leq k\leq p_n \}}- {\bf 1}_{\{\tau >k ; n-p_n-2\leq k\leq n-2 \}}\right]\\
&=&2^{n-1} (\mu_n\otimes{\cal B}_{n-1,1/2})\left[{\bf 1}_{\{\tau >k ; 0\leq k\leq n-3-p_n \}}- {\bf 1}_{\{\tau >k ; p_n+1\leq k\leq n-2 \}}\right]\\
&\sim& {\cal B}_{n-1,1/2}(0\leq k\leq n-3-p_n),
\end{eqnarray*}
where the equivalent is proved just as above.\\

\noindent
Finally, consider the case $p_n=n/2+\alpha_n\sqrt n$ with $\alpha_n\to\alpha$. We use the central limit Theorem and write
\begin{eqnarray*}
& &\mu_n\left[b_{2p_n}(N_\ell)\right]\\
&=&2^{n-1} \left(\mu_n\otimes{\cal B}_{n-1,1/2})\left[{\bf 1}_{\{\tau >k ;  k\leq p_n \}}- {\bf 1}_{\{\tau >k ; n-p_n-2\leq k\leq n-2 \}}\right]\right)\\
&=&2^{n-1} \left(\mu_n\otimes{\cal B}_{n-1,1/2})\left[{\bf 1}_{\{n^{-1/2}(\tau-n/2) >n^{-1/2}(k-n/2) ; n^{-1/2}(k-n/2)\leq \alpha_n \}} \right.\right.\\
& &\left.\left.-{\bf 1}_{\{n^{-1/2}(\tau-n/2) >n^{-1/2}(k-n/2) ; -\alpha_n-2n^{-1/2}\leq n^{-1/2}(k-n/2)\leq n^{-1/2}(n/2-2) \}}\right]\right)\\
&\sim& 2^{n-1} \mathbb E[{\bf 1}_{\{\sigma_\tau G_1 > G_2/2 ; G_2/2< \alpha \}}- {\bf 1}_{\{\sigma_\tau G_1 >G_2/2 ; G_2/2>-\alpha \}}]
\end{eqnarray*}
with $G_1$ and $G_2$ independent standard Gaussian random variables. The constant $C(\alpha)$ corresponds to the expectation in the last line. Using symetry properties for the distribution of $(G_1,G_2)$, we easily verify the announced formula for $C(\alpha)$. This ends the proof of Proposition \ref{prop:pr1-spatial}\CQFD

\section{Asymptotic behavior of the Poincar\'e polynomial\label{AP}}
\subsection{Planar polygons\label{APlanar}}
We will here consider the random Poincar\'e polynomial $p_{M_\ell}(t)$
as given in (\ref{PC2}) in the large $n$ limit. We first give a representation of this invariant in terms of random permutations and stopping times.

\begin{lemme}\label{lem:Poincaré} For any $t> 0$, the Poincaré polynomial is given $\mu_n$-almost surely by
$$p_{M_\ell}(t)=(1+t)^{n-1}({\cal U}_{\Sigma_{n-1}}\otimes {\cal B}_{n-1,\frac{t}{1+t}}) \left[  \ind_{\{\tau_\sigma(\tilde \ell) > k\}}+t^{-2}\ind_{\{ \tau_\sigma(\tilde\ell) > n-1-k\}}\right].$$
As a consequence,
$$\bar p_{M_\ell}(t)=(1+t)^{n-1}(\mu_n\otimes {\cal B}_{n-1,\frac{t}{1+t}}) \left[  \ind_{\{\tilde\tau > k\}}+t^{-2}\ind_{\{ \tilde\tau > n-1-k\}}\right].$$
\end{lemme}

Thanks to this lemma, we prove the following Proposition giving the asymptotic of the average Poincar\'e polynomial.
\begin{prop}\label{prop:pr2}
Let $\bar p_{M_\ell}(t)$ be the mean Poincar\'e polynomial. When $t>0$,
\begin{enumerate}
\item If $0<t<1$, then $\bar p_{M_\ell}(t)\sim (1+t)^{n-1}$.
\item If $t>1$, then $\bar p_{M_\ell}(t)\sim (1+t)^{n-1}t^{-2}$.
\item If $t=1$, then the mean total Betti number satisfies $\bar p_{M_\ell}(1)\sim 2^{n-1}.$
\end{enumerate}
\end{prop}

\noindent {\bf Proof of Lemma \ref{lem:Poincaré}} \ \\
Equation (\ref{PC2}) together with Lemma \ref{lem:representation1} yield
\begin{eqnarray*}
q(t) &=& \sum_{k=0}^{n-1} {{n-1}\choose k} t^k  {\cal U}_{\Sigma_{n-1}}(\tau_\sigma(\tilde \ell) >k) \\
&=& (1+t)^{n-1}\sum_{k=0}^{n-1} {{n-1}\choose k} \left(\frac{t}{1+t}\right)^k\left(\frac{1}{1+t}\right)^{n-1-k}  {\cal U}_{\Sigma_{n-1}}(\tau_\sigma(\tilde \ell) >k)\\
&=& (1+t)^{n-1}\left[({\cal U}_{\Sigma_{n-1}} \otimes {\cal B}_{n-1,\frac{t}{1+t}}) \left( \tau_\sigma(\tilde\ell) > k\right)\right].
\end{eqnarray*}
Please note that in the sum the terms corresponding to $k=n-2$ and $k=n-1$ vanish.
Finally,  Lemma \ref{lem:Poincaré} follows from the relation
$$ p_{M_\ell}(t)= q(t)+t^{n-3} q(t^{-1})+r(t),$$
with $r(t)$ $\mu_n$-almost surely vanishing and from the fact that the distribution of $k$ under ${\cal B}_{n-1,\frac{1}{1+t}}$ is equal to the distribution of $n-1-k$ under ${\cal B}_{n-1,\frac{t}{1+t}}$.\\
We use once again the invariance property of $\mu_n$ under the action of the symetric group to simplify the expression of the average Poincar\'e polynomial $\mu_n[p_{M_\ell}(t)]$. 
\CQFD\\
\ \\

\noindent
{\bf Proof of Proposition \ref{prop:pr2}}\\
We use the representation of the average Poincaré polynomial given in Lemma \ref{lem:Poincaré} together with weak convergence 
for $(\tilde\tau, k)$ under $\mu_n\otimes{\cal B}_{n-1,\frac{t}{1+t}}$ to study the asymptotic behavior .\\

\noindent
The weak law of large number for $\tilde\tau$ (see Lemma \ref{lem:LLN-TCL}) and a standard weak law of large numbers for binomial distribution imply that $(n^{-1}\tilde\tau, n^{-1}k)$ converges weakly under $\mu_n\otimes{\cal B}_{n-1,\frac{t}{1+t}}$
to $(0,\frac{t}{1+t})$. The continuous mapping theorem implies that for $0<t<1$ or $t>1$, the following weak convergence holds
under $\mu_n\otimes{\cal B}_{n-1,\frac{t}{1+t}}$:
$$
\ind_{\{\tilde\tau > k\}}+t^{-2}\ind_{\{ \tilde\tau > n-1-k\}}\Rightarrow \ind_{\{\frac{1}{2} > \frac{t}{1+t}\}}+t^{-2}\ind_{\{ \frac{1}{2} > 1-\frac{t}{1+t}\}}=\min(1,t^{-2}).
$$
Integrating this (bounded) convergence yield the result for $t\neq 1$.\\

\noindent
For $t=1$, the continuous mapping theorem does not hold no longer since the map $(\tilde\tau,k)\mapsto \ind_{\{\tilde\tau > k\}}$ is not continuous at point $(1/2,1/2)$. We need here the central limit Theorem. From Lemma \ref{lem:LLN-TCL} and standard results for binomial distribution, $(n^{-1/2}(\tilde\tau-n/2), n^{-1/2}(k-n/2))$ converges weakly under $\mu_n\otimes{\cal B}_{n-1,1/2}$ to ${\cal N}(0,\sigma_{\tilde\tau}^2)\otimes{\cal N}(0,1/4)$. The continuous mapping Theorem yields
\begin{eqnarray*}
& &\ind_{\{\tilde\tau > k\}}+\ind_{\{ \tilde\tau > n-1-k\}}\\
&=&\ind_{\{n^{-1/2}(\tilde\tau-n/2)>n^{-1/2}(k-n/2)\}}+\ind_{\{ n^{-1/2}(\tilde\tau-n/2), n^{-1/2}(n/2-1-k)\}}\\
&\Rightarrow& \ind_{\{\sigma_{\tilde\tau}G_1>G_2/2\}}+\ind_{\{\sigma_{\tilde\tau}G_1>-G_2/2\}}
\end{eqnarray*}
with $G_1$ and $G_2$ independent standard Gaussian random variables. We integrate this (bounded) convergence and remark
that ${\mathbb E}(\ind_{\{\sigma_{\tilde\tau}G_1>G_2/2\}})={\mathbb E}(\ind_{\{\sigma_{\tilde\tau}G_1>-G_2/2\}})=1/2$.
\CQFD \\
\ \\

\noindent
{\bf Remark:} we can use large deviations results to estimate the speed of convergence in Proposition \ref{prop:pr2} when $t\neq 1$. For example for $0<t<1$, write
\begin{eqnarray*}
& &\mu_n \left[ (1+t)^{-(n-1)}p_{M_\ell}(t)-1 \right]  \\
&=& (\mu_n\otimes {\cal B}_{n-1,\frac{t}{1+t}})\left[ {\bf 1}_{\{\tilde\tau >k \}}- 1+t^{-2}{\bf 1}_{\{\tilde\tau >n-1-k \}}\right]\\
&=& (\mu_n\otimes {\cal B}_{n-1,\frac{t}{1+t}})\left[ {\bf 1}_{\{n^{-1}(\tilde\tau-k) \leq 0 \}}+t^{-2}{\bf 1}_{\{n^{-1}(\tilde\tau+k) \geq 1 \}}\right].
\end{eqnarray*}
Now large deviations  for $n^{-1}(\tilde \tau,k)$ under $(\mu_n\otimes {\cal B}_{n-1,\frac{t}{1+t}})$
will give the speed of convergence to $0$ in a logarithmic scale.\\
For $t>1$, we have
\begin{eqnarray*}
& &\mu_n \left[ (1+t)^{-(n-1)}(p_{M_\ell}(t))-t^{-2} \right]  \\
&=& (\mu_n\otimes {\cal B}_{n-1,\frac{t}{1+t}})\left[ {\bf 1}_{\{n^{-1}(\tilde\tau-k) > 0 \}}+t^{-2}{\bf 1}_{\{n^{-1}(\tilde\tau+k) < 1 \}}\right],
\end{eqnarray*}
and we can use the same method.

\subsection{Spatial polygons\label{Spa}}

We use the same strategy in the spatial case and use formula (\ref{PC3}) giving the Poincaré polynomial for generic vector length. Since $\mu$ is diffuse, $\mu_n$-almost every vector length is generic and equation (\ref{PC3}) holds.
The related total Betti number is obtained by taking the $t\to 1$ limit in (\ref{PC3})
\begin{eqnarray}
p_{N_\ell}(1)&=&\lim_{t\to 1}\frac{1}{1-t^2}\Big(\hat q(t^2)-t^{2(n-2)}\hat q(t^{-2})\Big)\nonumber \\
           &=&(n-2)\hat q(1)-2\hat q'(1)\nonumber \\
	&=& (n-2)\sum_{j=0}^{n-1} \hat a_j -2\sum_{j=0}^{n-1}j\hat a_j \label{TotalBetti}
\end{eqnarray}

We use the following representations for the Poincaré polynomial:
\begin{lemme}\label{lem:Poincaré2}
The Poincaré polynomial is given $\mu_n$-almost surely by
$$
p_{N_\ell}(t) =\frac{(1+t^2)^{n-1}}{1-t^2}( {\cal U}_{\Sigma_{n-1}}\otimes {\cal B}_{n-1,\frac{t^2}{1+t^2}})\left[ {\bf 1}_{\{\tau_\sigma(\ell)>k\}}-t^{-2}{\bf 1}_{\{\tau_\sigma(\ell)>n-k\}}\right], 
$$
for $0<t<1$ or $t>1$, and by
$$
p_{N_\ell}(1) =n2^{n-1}({\cal U}_{\Sigma_{n-1}} \otimes {\cal B}_{n-1,1/2}) \left[ (\frac{n-2}{n}-\frac{2k}{n}){\bf 1}_{\{\tau_\sigma(\ell) > k\}}\right]$$
for $t=1$.

\end{lemme}

\begin{prop}\label{prop:pr3}
Let $\bar p_{N_\ell}$ be the mean Poincar\'e polynomial associated with random spatial polygons. When $t>0$,
\begin{enumerate}
\item If $0<t<1$, then $\bar p_{N_\ell}(t)\sim \frac{(1+t^2)^{n-1}}{1-t^2}$.
\item If $t>1$, then $\bar p_{N_\ell}(t)\sim \frac{(1+t^2)^{n-1}}{t^2(t^2-1)}$.
\item If $t=1$, then the total Betti number satisfies $\bar p_{N_\ell}(1)\sim  n 2^{n-2}.$
\end{enumerate}
\end{prop}
\medskip

\noindent {\bf Remark:} In the case of spatial polygons, the Poincaré polynomial is an even function. Hence its asymptotic mean behavior for $t<0$ follows directly from Proposition \ref{prop:pr3}.\\

\noindent {\bf Proof of Lemma \ref{lem:Poincaré2}}\\
The proof is very similar to the proof of Lemma \ref{lem:Poincaré}. Equation (\ref{PC3}) together with Lemma \ref{lem:representation2} yield
\begin{eqnarray*}
\hat q(t) &=& (1+t)^{n-1}({\cal U}_{\Sigma_{n-1}} \otimes {\cal B}_{n-1,\frac{t}{1+t}}) \left[ \tau_\sigma(\ell) > k\right].
\end{eqnarray*}
The case $t\neq 1$ follows from the relation 
$$ p_{N_\ell}(t)=\frac{1}{1-t^2}\Big(\hat q(t^2)-t^{2(n-2)}\hat q(t^{-2})\Big)$$
and from the fact that the distribution of $k$ under ${\cal B}_{n-1,\frac{1}{1+t^2}}$ is equal to the distribution of $n-K$ under 
${\cal B}_{n-1,\frac{t^2}{1+t^2}}$.\\

\noindent
In the case $t=1$, equation (\ref{TotalBetti}) and Lemma \ref{lem:representation2} imply
\begin{eqnarray*}
p_{N_\ell}(1)&=&(n-2)\sum_{j=0}^{n-1} \hat a_j -2\sum_{j=0}^{n-1}j\hat a_j \\
&=& n2^{n-1}\sum_{j=0}^{n-1} (\frac{n-2}{n}-\frac{2j}{n}){ n-1\choose j} {\cal U}_{\Sigma_{n-1}}[\tau_\sigma(\ell)>j]\\
&=& n2^{n-1}({\cal U}_{\Sigma_{n-1}} \otimes {\cal B}_{n-1,1/2}) \left[ (\frac{n-2}{n}-\frac{2k}{n}){\bf 1}_{\{\tau_\sigma(\ell) > k\}}\right]
\end{eqnarray*}
\CQFD
\ \\

\noindent
{\bf Proof of Proposition \ref{prop:pr3}}
The case $0<t<1$ and $t>1$ are easily deduced from Lemma \ref{lem:Poincaré2} using the following law of large numbers: under ${\cal B}_{n-1,\frac{t^2}{1+t^2}}\otimes {\cal U}_{\Sigma_{n-1}}$, $n^{-1}(k,\tau)$ converges weakly to $(1/2,\frac{t^2}{1+t^2})$ as $n\to\infty$.
Details are omitted since they are as in the proof of Proposition \ref{prop:pr2}.\\

In the case $t=1$, the central limit Theorem from Lemma \ref{lem:LLN-TCL} states that $(n^{-1/2}(\tilde\tau-n/2), n^{-1/2}(k-n/2))$ converges weakly under $\mu_n\otimes{\cal B}_{n-1,1/2}$ to ${\cal N}(0,\sigma_{\tilde\tau}^2)\otimes{\cal N}(0,1/4)$. As a consequence,
\begin{eqnarray*}
n^{-1}2^{-(n-1)} \mu_n[p_{N_\ell}(1)] &=&(\mu_n \otimes {\cal B}_{n-1,1/2}) \left[ (\frac{n-2}{n}-\frac{2k}{n}){\bf 1}_{\{\tau > k\}}\right]\\
&\rightarrow & {\mathbb E}\left[1- 2\frac{1}{2}{\bf 1}_{\{\sigma_\tau G_1 > G_2/2\}}\right]=\frac{1}{2},
\end{eqnarray*}
whith $G_1$ and $G_2$ independent standard Gaussian random variables.\CQFD

\ \\

\noindent
{\bf Remark:} In order to estimate the speed of convergence in Proposition {\ref{prop:pr3} when $t\neq 1$, we can use for $0<t<1$ the expression
\begin{eqnarray*}
& &\mu_n \left[ (1-t^2)(1+t^2)^{-(n-1)}p_{M_\ell}(t)-1 \right]\\
&=& (\mu_n\otimes {\cal B}_{n-1,\frac{t^2}{1+t^2}})\left[ {\bf 1}_{\{ n^{-1}(\tau-k)\leq 0\}}-t^{-2}{\bf 1}_{\{ n^{-1}(\tau+k)> 1\}}\right]
\end{eqnarray*}
and  for $t>1$ 
\begin{eqnarray*}
& &\mu_n \left[ (t^2-1)(1+t^2)^{-(n-1)}p_{M_\ell}(t)-t^{-2} \right]\\
&=& (\mu_n\otimes {\cal B}_{n-1,\frac{t^2}{1+t^2}})\left[ {\bf 1}_{\{ n^{-1}(\tau-k) > 0\}}-t^{-2}{\bf 1}_{\{ n^{-1}(\tau+k)\leq 1\}}\right]. 
\end{eqnarray*}

Large deviations results for $n^{-1}(\tau,k)$ under $(\mu_n\otimes {\cal B}_{n-1,\frac{t^2}{1+t^2}})$ would give the speed of convergence in a logarithmic scale.

\subsection{Higher moments}
We consider here the higher moments of the Poincaré polynomial and prove that their asymptotic 
behavior is given by the first moment. To this aim, we prove a weak law of large numbers for the renormalized Poincaré polynomial.\\

\noindent
We begin with the case of planar polygon.

\begin{prop}\label{prop:pr4}
For any $t>0$, the following weak convergence holds under $\mu_n$ as $n\to\infty$
$$(1+t)^{-(n-1)}p_{M_\ell(t)} \Rightarrow \min(1,t^{-2}).$$
As a consequence, for any $t>0$ and $\nu\in{\mathbb N}$,
$$\mu_n\left[ p_{M_\ell}(t)^\nu\right] \sim \left(\mu_n\left[ p_{M_\ell}(t)\right] \right)^\nu$$
\end{prop}
{\bf Proof of Proposition \ref{prop:pr4}}\\
Proposition \ref{prop:pr2} states that the expectation under $\mu_n$ of $(1+t)^{-(n-1)}p_{M_\ell(t)}$ converges to $\min(1,t^{-2})$ as $n\to\infty$. Hence, weak convergence will be proved as soon as we show that the variance under $\mu_n$ of
$(1+t)^{-(n-1)}p_{M_\ell(t)}$ goes to zero. We use the representation of the Poincaré polynomial from Lemma \ref{lem:Poincaré}
and the replica trick to compute the second moment
$$\mu_n\left[(1+t)^{-2(n-1)}p_{M_\ell}(t)^2\right]=(\mu_n\otimes {\cal B}_{n-1,\frac{t}{1+t}}^{\otimes 2}\otimes {\cal U}_{\Sigma_{n-1 }}^{\otimes 2})\Big[{\rm Prod} \Big].$$
with
$${\rm Prod}=({\bf 1}_{\{\tau_{\sigma_1}(\tilde\ell)>k_1\}} +t^{-2} {\bf 1}_{\{\tau_{\sigma_1}(\tilde\ell)>n-1-k_1\}})
({\bf 1}_{\{\tau_{\sigma_2}(\tilde\ell)>k_2\}}+t^{-2} {\bf 1}_{\{\tau_{\sigma_2}(\tilde\ell)>n-1-k_2\}}). $$
We need to show that the two factors of ${\rm Prod}$ are asymptotically independent in the limit $n\to\infty$. This would yield
$$
\mu_n\left[(1+t)^{-2(n-1)}p_{M_\ell}(t)^2\right]\sim \left(\mu_n\left[(1+t)^{-(n-1)}p_{M_\ell}(t)\right]\right)^2,
$$
and hence the variance of $\frac{p_{M_\ell}(t)}{(1+t)^{n-1}}$ would converge to zero as $n\to\infty$. 
We now prove asymptotic independance of the two factors. 
When $0<t<1$ or $t>1$ the aymptotic independence follows from the weak law of large numbers obtained in  Lemma \ref{lem:LLN-TCL}, both factors converging weakly to $\min(1,t^{-2})$ (note that the distribution of $\tau_{\sigma}(\tilde\ell)$ under $\mu_n\otimes{\cal U}_{\Sigma_{n-1 }}$ is equal to the distribution of $\tilde\tau(\ell)$ under $\mu_n$). When $t=1$, we use the bivariate central limit Theorem stated in Lemma \ref{lem:TCL2} in the Appendix.
Weak convergence is proved.\\

The convergence of the moments is a direct consequence of the weak convergence once we remark that the renormalized Poincaré polynomial $(1+t)^{-(n-1)}p_{M_\ell}(t)$ is $\mu_n$ almost surely bounded by $1+t^{-2}$ (this is clear from the representation given in Lemma \ref{lem:Poincaré}).\CQFD\\
\ \\

\vspace{1cm}

We consider now the higher moments of the Poincaré polynomial for spatial polygons spaces. The results and methods are very similar to one of the planar case and are based on Lemma \ref{lem:Poincaré2}. Hence we give only the main lines of the proof.
\begin{prop}\label{prop:pr5}
The following weak convergence holds under $\mu_n$ as $n\to\infty$,
\begin{eqnarray*}
{\rm if}& 0<t<1,\ \ & (1+t^2)^{-(n-1)}p_{N_\ell}(t)\Rightarrow (1-t^2)^{-1},\\  
{\rm if}& t>1,\ \ & (1+t^2)^{-(n-1)}p_{N_\ell}(t)\Rightarrow t^{-2}(t^2-1)^{-1},\\
{\rm if}& t=1, \ \ & n^{-1}2^{-n}p_{N_\ell}(1)\Rightarrow 1/4.
\end{eqnarray*}
As a consequence, for any $t>0$ and $\nu\in{\mathbb N}$,
$$\mu_n\left[ p_{N_\ell}(t)^\nu\right] \sim \left(\mu_n\left[ p_{N_\ell}(t)\right] \right)^\nu.$$
\end{prop}

\noindent
{\bf Proof of Proposition \ref{prop:pr5}}\\
The proof is similar to the proof of Proposition \ref{prop:pr4} with the following expression 
of the renormalized Poincaré polynomial deduced from Lemma \ref{lem:representation2}: for $0<t<1$ or $t>1$
$$
(1-t^2)(1+t^2)^{-(n-1)}p_{N_\ell}(t) =( {\cal U}_{\Sigma_{n-1}}\otimes {\cal B}_{n-1,\frac{t^2}{1+t^2}})\left[ {\bf 1}_{\{\tau_\sigma(\ell)>k\}}-t^{-2}{\bf 1}_{\{\tau_\sigma(\ell)>n-k\}}\right], 
$$
and for $t=1$
$$
n^{-1}2^{-(n-1)}p_{N_\ell}(1) =({\cal U}_{\Sigma_{n-1}} \otimes {\cal B}_{n-1,1/2}) \left[ (\frac{n-2}{n}-\frac{2k}{n}){\bf 1}_{\{\tau_\sigma(\ell) > k\}}\right]
$$
Convergence of the expectation was proved in Proposition \ref{prop:pr3}. The variance is computed using thanks to the replica trick and is shown to converge to zero because of the asymptotic independence of 
${\bf 1}_{\{\tau_{\sigma_i}(\ell)>k_i\}}, i=1,2$ under $\mu_n\otimes {\cal B}_{n-1,\frac{t^2}{1+t^2}}^{\otimes 2}\otimes {\cal U}_{\Sigma_{n-1}}^{\otimes 2}$ (see Lemma\ref{lem:TCL2}). \CQFD

\section*{Appendix}
\noindent {\bf Proof of Lemma \ref{lem:LLN-TCL}}\\
The weak law of large number is a consequence of the central limit theorem that we prove now.
Let $p_n=\frac{n}{2}+\alpha_n\sqrt n$ with $\alpha_n\to\alpha$ as $n\to\infty$. Using the definition of $\tilde\tau$, 
\begin{eqnarray*}
\mu_n\left(\tilde\tau\leq p_n\right)
&=&\mu_n\left(\tilde l_n+\sum_{i=1}^{p_n} \tilde l_i-\sum_{i=p_n+1}^{n-1} \tilde l_i\geq 0\right)\\
&=& \mu_n\left(n^{-1/2}\tilde l_n+n^{-1/2}(\sum_{i=1}^{p_n}  \tilde l_i-\sum_{i=p_n+1}^{n-1} \tilde l_i)\geq 0 \right).
\end{eqnarray*}
We now prove that $n^{-1/2}\tilde l_n$ converges weakly to zero and that $n^{-1/2}(\sum_{i=1}^{p_n}  \tilde l_i-\sum_{i=p_n+1}^{n-1} \tilde l_i)$ satisfies a central limit theorem.
To see this, we denote by $F_\mu$ the repartition function of $\mu$, and remark that the distribution of $\tilde l_n$ is given by
$$
\mu_n(\tilde l_n\leq x)=F_\mu(x)^n.
$$
Hence 
$$\mu_n(n^{-1/2}\tilde l_n > \varepsilon )=(1-F_\mu(\varepsilon n^{1/2}))^n,$$
and the exponential Markov inequality implies
$$1-F_\mu(\varepsilon n^{1/2})\leq \exp(-\eta \varepsilon n^{1/2}) \int e^{\eta x}\mu(dx),\ \eta > 0.$$
This implies the weak convergence $n^{-1/2}\tilde l_n$ to zero.
Conditionnaly to $\tilde l_n=u$, the other components $(l_i)_{1\leq i\leq n-1}$ are i.i.d. with conditional distribution given by
$$\mu_n(\tilde l_i \leq x \ |\ \tilde l_n=u)= \frac{F_\mu(x\wedge u)}{F_\mu(u)}.$$
Denote by $m(u)$ and $\sigma^2(u)$ the related conditionnal expectation and variance.
From the central limit theorem for independent variables, conditionnaly to $\tilde l_n=u$, the quantity 
$n^{-1/2}(\sum_{i=1}^{p_n}  \tilde l_i-\sum_{i=p_n+1}^{n-1} \tilde l_i)$ converges weakly to a gaussian distribution of mean $2\alpha m(u)$ and variance $\sigma^2(u)$. Hence the conditionnal probability
$$\mu_n\left[n^{-1/2}(\sum_{i=1}^{p_n}  \tilde l_i-\sum_{i=p_n+1}^{n-1} \tilde l_i)\geq 0  \vert \tilde l_n=u\right] $$
converges to $F_{\cal N}(2\alpha m(u)/\sigma(u))$ as $n\to \infty$. We now have to integrate this with respect to $\tilde l_n$. Taking into account that $\tilde l_n$ converges weakly to $l_{\rm max}=\inf \{x\in\real;F_\mu(x)=1\}\in (0,+\infty]$ as $n\to\infty$ and that $(m(u),\sigma(u))\to (m,\sigma)$ as $u\to  l_{\rm max}$, we see that 
$$\mu_n\left[n^{-1/2}(\sum_{i=1}^{p_n}  \tilde l_i-\sum_{i=p_n+1}^{n-1} \tilde l_i)\geq 0  \right]\to F_{\cal N}(2\alpha m/\sigma). $$ 
This proves the central limit theorem for $\tilde\tau$.  

We now prove the large deviation estimate. Since 
$$\mu_n\left(\tilde\tau\leq (1/2-\varepsilon) n)\right)=\mu_n\left(\tilde l_n+\sum_{i=1}^{[(1/2-\varepsilon)n]} \tilde l_i-\sum_{i=[(1/2-\varepsilon)n]}^{n-1} \tilde l_i\geq 0\right),$$
we will provide large deviations estimates for the random sum
$$S_n=\tilde l_n+\sum_{i=1}^{[(1/2-\varepsilon)n]} \tilde l_i-\sum_{i=[(1/2+\varepsilon)n]}^{n-1} \tilde l_i.$$
For $t\in\mathbb R$, the logarithmic moment generating function is defined by
$$
\Lambda_n(t)=\log(\mu_n(\exp(tS_n))).
$$  
Using Laplace method, we see that as $n\to\infty$, $n^{-1}\Lambda_n(t)$ converges to 
$$
\Lambda(t)=(1/2-\varepsilon)\int e^{ty}\mu({\rm d}y)+(1/2+\varepsilon)\int e^{ty}\mu({\rm d}y).
$$
Using Gärtner-Ellis theorem, see e.g. \cite{Dembo}, we deduce a large deviations principle for the sum $n^{-1}S_n$ of speed $n$ and of good rate function $I$ being the Fenchel-Legendre transform of $\Lambda$. The exact form of $I$ is irrelevant here but it is important to see that $I$ is strictly positive on $[0,\infty)$. Standard arguments from large deviations theory (see \cite{Dembo}) give that $I$ vanishes only at $(1/2-\varepsilon)m-(1/2+\varepsilon)m=-2\varepsilon m<0$, and hence the action $I$ is negative on $[0,\infty)$. As a consequence, the large deviations principle states that
$$\limsup n^{-1}\log\mu_n\left(\tilde\tau\leq (1/2-\varepsilon) n\right)\leq -\inf_{[0,\infty)} I <0.$$
The same technique is used to deal with $\mu_n\left(\tilde\tau\geq (1/2+\varepsilon) n\right)$ and this proves the Lemma.
\CQFD

\begin{lemme}\label{lem:TCL2}
The following bivariate Central Limit Theorem holds under $\mu_n\otimes {\cal U}_{\Sigma_{n-1 }}^{\otimes 2}$:
$$n^{-1/2}( \tau_{\sigma_1}(\ell)-n/2,\tau_{\sigma_2}(\ell)-n/2) \Rightarrow {\cal N}(0,\sigma_\tau^2)^{\otimes 2}.$$
It also holds for $\tilde\tau$

\end{lemme}
{\it Proof of Lemma }
need a bivariate central limit Theorem for  $(\tau_{\sigma_1}(\tilde\ell),\tau_{\sigma_2}(\tilde\ell))$ under  
Let $p_{n,i}=\frac{n}{2}+\alpha_{n,i}\sqrt n$ with $\alpha_{n,i}\to\alpha$ as $n\to\infty$ for $i=1,2$. By the definition of $\tilde\tau_{\sigma_i}$, 
\begin{eqnarray*}
& &(\mu_n\otimes {\cal U}_{\Sigma_{n-1 }}^{\otimes 2})\left(\tilde\tau_{\sigma_i}\leq p_{n,i} ; i=1..2\right)\\
&=& (\mu_n\otimes {\cal U}_{\Sigma_{n-1 }}^{\otimes 2})\left(n^{-1/2}\tilde l_n+n^{-1/2}(\sum_{j=1}^{p_{n,i}}  \tilde l_{\sigma_i(j)}-\sum_{j=p_{n,i}+1}^{n-1} \tilde l_{\sigma(j)})\geq 0 ; i=1,2\right).\\
\end{eqnarray*}
We know from the proof of Lemma \ref{lem:LLN-TCL} that  $n^{-1/2}\tilde l_n$ converges weakly to zero. It remains to check that $n^{-1/2}(\sum_{j=1}^{p_{n,i}}  \tilde l_{\sigma_i(j)}-\sum_{j=p_{n,i}+1}^{n-1} \tilde l_{\sigma_i(j)})_{i=1,2}$ satisfies a bivariate central limit theorem. Let $\theta_i, i=1,2$ be real numbers, and consider the linear combination
$$
\sum_{i=1}^2 \theta_i n^{-1/2}(\sum_{j=1}^{p_{n,i}}  \tilde l_{\sigma_i(j)}-\sum_{j=p_{n,i}+1}^{n-1} \tilde l_{\sigma_i(j)}
 = n^{-1/2}\sum_{j=1}^{n-1} (\theta_1\varepsilon_{n,1}(j)+\theta_1\varepsilon_{n,2}(j))\tilde l_j,
$$
where we set
$\varepsilon_{n,i}(j)=2{\bf 1}_{\{\sigma_i(j)\leq p_{n,i}\}}-1$.
Conditionnaly to $\tilde l_n=u$, the components $\tilde l_j$ are i.i.d. with mean $m(u)$ and variance $\sigma(u)$, and hence the above sum is a linear triangular array of independent variables with random coefficients $(\theta_1\varepsilon_{n,1}(j)+\theta_1\varepsilon_{n,2}(j))_{1\leq j\leq n-1}$. The coefficients are almost surely bounded and satisfy a weak law of large numbers under ${\cal U}_{\Sigma_{n-1 }}^{\otimes 2}$
$$n^{-1}\sum_{j=1}^{n-1}(\theta_1\varepsilon_{n,1}(j)+\theta_1\varepsilon_{n,2}(j))^2\to \theta_1^2+\theta_2^2.$$ 
(note that the empirical distribution $\frac{1}{n-1}\sum_{j=1}^{n-1} \delta_{(\varepsilon_{n,1}(j),\varepsilon_{n,1}(j)\varepsilon_{n,2}(j))}$ converges weakly to the uniform distribution on $\{(\pm 1,\pm 1)\}$). As a consequence, conditionaly to $\tilde l_n=u$, the above sum converges to a gaussian random variables of mean $2(\alpha_1\theta_1+\alpha_2\theta_2)m(u)$ and variance $(\theta_1^2+\theta_2^2)\sigma^2(u)$. Integrating with respect to $\tilde l_n$ we obtain that the sum converges weakly to a gaussian random variables with mean $2(\alpha_1\theta_1+\alpha_2\theta_2)m$ and variance $(\theta_1^2+\theta_2^2)\sigma^2$. This proves the bivariate central limit theorem with asymptotic independent components.
\CQFD

\end{document}